\magnification=1200
\noindent


\bigskip

\centerline{\bf Arithmetic of rationaly connected quintic $3$-folds}
\centerline{\bf  over finite and functions fields}
\bigskip

\centerline{Marc Perret}

\medskip


\bigskip

 \noindent {\bf I. Introduction.} The class of Rationaly Connected (RC for short) algebraic varieties was introduced independantly by F. Campana in [2] and
 by Kollar, Miyahoka and Mori in [10] in 1992.
 This is a natural extension in dimension greater than 2 of the class
 of rational curves. Roughtly speaking,
 a variety $X$ defined over a field $k$ is rationaly connected if there exists some dense open subset $U$ of $X$, such that for any $P, Q \in U$ there exists $\varphi : {\bf P}^1_k \rightarrow X$, such that $P=\varphi(0)$ and $Q=\varphi(\infty)$. Formally, we say that $X$ is rationally connected if there exists a familly of proper curves
 $g : U \rightarrow Y$ whose geometric fibers are irreducible rational curves, with a cycle morphism
 $u : U \rightarrow X$ such that $u^{(2)} : U\times_YU \rightarrow X \times X$ is dominant
 (see [9]).

 \medskip
 
Several families of RC variety are known. For $n \geq 2$ and $d \geq 1$ two integers, we denote by $X_d^{(n-1)}$ an hypersurface of ${\bf P}^n_k$ of degree $d$, defined over $k$.

 \medskip
\noindent {\it Type I.} Rational varieties, and more generally unirational varieties, are RC.
  
  \medskip
  
 \noindent{\it Type II.} $X_d^{(n-1)} \subset {\bf P}^n$ is RC if $d \leq n$.
 
 \medskip 
 
\noindent {\it Type III.} It is well known (see for instance [12]) that if  $k$ is algebricaly closed of characteristic prime to  $d$, then 
 $$\pi_1({\bf P}^n-X_d^{(n-1)}) \simeq {\bf Z}/d{\bf Z}.$$
 If $d'$ divides $d$, there is an unique morphism $Y_{X, d'} \rightarrow {\bf P}^n$ of degree $d'$ of ${\bf P}^n$, ramified along $X$. If $k$ is not algebraically closed, $Y_{X, d'}$ possesses several forms on $k$. We have then :
 
$$Y_{X, d'} \rm{~is~RC~if~} n+d'\leq d.$$
 
 \medskip
 
\noindent {\it Type IV.} Let $X=X_d^{(n-1)} \subset {\bf P}^n$. Suppose that $X$ contains a line $L=X_1^{(1)} \subset X \subset {\bf P}^n$, with multiplicity $d-2$. Denote by $\widetilde X_d^{(n-1)}$ the strict transform of $X_d^{(n-1)}$ in the blowing up of ${\bf P}^n$ along $L$. Then
 $$\widetilde X_d^{(n-1)} \rm{~is~RC}.$$

 \bigskip
 
 If $X$ is a variety defined over a field $k$, we denote by $X(k)$ the set of rational points of $X$ over $k$, and by $\sharp X(k)$ its cardinal. We are specially interested with both cases where  $k$ is finite in chapter III, and where $k$ is the rational function field over an algebraically closed field in chapter IV. A central problem of the arithmetic of varieties is:
 
 \medskip
 
\noindent {\bf Problem.} {\it Is $X(k)$ non-empty ?}

 \bigskip
 
 Consider first the case of hypersurfaces of ${\bf P}^n$. Then the three following theorems are well known:
 
 \medskip
 
\noindent  {\bf Theorem (Chevalley-Warning, 1936)} {\it let $P(X_0, \dots, X_n)$ be a polynomial of degree $d$, with coefficients in a finite field $k$ of characteristic $p >0$. Suppose that $d \leq n$. Then the number of solutions $(x_0, \dots, x_n)$ in $k^{n+1}$ of the equation $P(x_1, \dots, x_n)=0$ is a multiple of $p$.}
 
 \bigskip
 
\noindent  {\bf Corollary.} {\it In the situation of Chevalley-Warning theorem, if we suppose moreover that  $P$ is homogeneous, then the number of  $k$-rational points of the projective variety $X$ defined by $P$ satisfies
 $$\sharp X(k) \equiv 1~(mod~p).$$}

 \medskip
 
This theorem was sharpened by Ax in [1] (and then by Katz):
 
 \medskip
 
\noindent  {\bf Theorem (Ax, 1964)}  {\it Let $P(X_0, \dots, X_n)$ be a polynomial of degree $d$, with coefficients in the finite field $k$ with $q$ elements. Then the number of  solutions $(x_0, \dots, x_n)$ in $k^{n+1}$ of the equation $P(x_1, \dots, x_n)=0$ is a multiple of $q^\alpha$, where 
 $$\alpha = \Big\lceilÊ{n-d \over d} \Big\rceil.$$}
 
 \bigskip
 
 Of course, Ax theorem admits a projective version if $P$ is homogeneous.
 For rational functions fields, we have:
  \medskip

\noindent  {\bf Theorem (Tsen, 1936)} {\it Let $k$ be an algebraically closed field, with rational function field  $k(T)$, and $P(X_0, \dots, X_n)$ homogeneous of degree $d < n$ with coefficients in $k(T)$. 
 Then there exists a non-zero solution $(f_0(T), \dots, f_n(T))$ of the equation
 $P(f_0(T), \dots, f_n(T)) = 0$.}

 \medskip

 Hence, if $d \leq n$, any $X^{(n-1)}_d \subset {\bf P}^n$ possess a rational point over a finite field, and any fibred variety in $X^{(n-1)}_d$ over the projective line over an algebraically closed field possess a section. Chevalley-Warning and Tsen theorems where recently generalised for smooth projective RC varieties. The proofs use deep notions, and are completely different in nature (see [3] for a proof of Esnault theorem, and [6] for a proof of Gabber and all theorem).

\bigskip

\noindent  {\bf Theorem (Esnault, 2002).} {\it Let $X$ be a proper smooth RC variety defined over the finite field ${\bf F}_q$ with $q$ elements. Then 
 $$\sharp X({\bf F}_q) \equiv 1~(mod~q).$$
In particular, $X({\bf F}_q)$ is non-empty. 
 }
 
 \bigskip
 
\noindent  {\bf Theorem (Gabber, Harris, Star, de Jong (2004))} {\it Let $X$ be a projective smooth RC variety over the field ${\bf C}(T)$. Then $X({\bf C}(T))$ is non-empty.}

  \bigskip
  
As already seen, up to the smoothness assumption,  Ax and Tsen theorems are particular cases of Esnault and Gabber and all ones for RC varieties of type II. In the same way, they are proved in [9] page 232 (for Tsen type theorem) and in [11] (for Ax type theorem) for varieties of type III, also without the smoothness hypothesis.
 
 The aim of the present paper is to prove Esnault and Gabber and all theorems for RC varieties of type IV, without the smoothness hypothesis. The nice points are that they give unified proofs, in the spirit of both original Ax's and Tsen's ones, and works in the singular case. The bad point is that it works only for sub-varieties of ``small degrees"  in an ad-hoc simplicial toric variety ${\bf P}_\Sigma$.
 
 For the sake of simplicity, we will focus on the first interesting case, namely the case of the blowing up of a quintic $3$-fold in ${\bf P}^4$ containing a line with multiplicity $3$, which justifies the title of this paper.

 \bigskip

 \noindent {\bf II. Toric varieties}
The Aim of this chapter is to introduce the example of section II. 3. 3, which will be crucial for our purpose. 
 Let $k$ be any field of characteristic zero or $p>0$. It will be a finite field in section III, and the field of rational functions over an algebraically closed field in section IV.

 \medskip
 
\noindent  {\bf II-1. Toric coordinate ring.} We will recall in this chapter the definition due to Cox, of simplicial toric varieties. The reader may refer to [4] for details.
 
To begin with, recall that the simplest example of proper toric variety is the projective space ${\bf P}^n$. Two different constructions may be given:
 
 - by gluing $n+1$ affine spaces ${\bf A}^n$ ;
 
- or  as a quotient of ${\bf A}^{n+1}-\{0\}$ by the multiplicative group ${\bf G}_m$.
 
 \medskip
 
 One can associate to a fan (see before) $\Sigma$ a toric variety ${\bf P}_\Sigma$, usualy by gluing affine spaces, see for instance [8] or [5].
If the fan $\Sigma$ is {\it simplicial}, one can alternatively define ${\bf P}_\Sigma$ as the quotient of an open subset of an affine space by a torus. This last point of vue of Cox enables us to speak about {\it homogeneous coordinates}, as in the projective space case.

 \medskip

\noindent  {\bf Definitions.} {\it 
 ($\imath$) A convex polyedral cone $\sigma$ in the vector space ${\bf R}^d$ is a subset of the form
 
 $$\sigma_{\{v_1, \dots, v_r\}} = \{\lambda_1 v_1 + \dots + \lambda_r v_r ; \lambda_1, \dots, \lambda_r \geq 0\}$$
where $v_1, \dots, v_r$ are $r$ vectors in ${\bf R}^d$. We say that $v_1, \dots, v_r$ generates the cone. 

 \noindent ($\imath\imath$) The cone $\sigma$ is simplicial if it is generated by a free subset in ${\bf R}^d$.
 
 \noindent ($\imath \imath \imath$) A fan $\Sigma$ in ${\bf R}^d$ is a finite set of polyedral convex cones, such that: 
  
 - Any face of a cone $\sigma \in \Sigma$ is also in $\Sigma$.
 
 - the intersection of two cones in $\Sigma$ is a face of them.

 \noindent ($\imath v$) The fan $\Sigma$ is simplicial if all its cones are.
  }
 
 \bigskip

\noindent  {\bf Definitions.} {\it
 $(\imath)$ Let $\Sigma$ be a simplicial fan in ${\bf R}^d$. Denote by $n_1, \dots, n_\rho$ the generators of its cones of dimension $1$. We define :

 $$G_\Sigma= \{(\mu_1, \dots, \mu_\rho) \in {\bf G}_m^\rho ;
 \prod_{i=1}^\rho \mu_i^{<n_i \vert x>}=1 \forall x \in {\bf R}^d\} \subset {\bf G}_m^\rho$$
 
 \noindent $(\imath \imath)$ A set $\{n_{{i_1}}, \dots, n_{{i_k}}\}$ is primitive if it is contained in no cone of  $\Sigma$, but each strict subsets are. If $\{n_{{i_1}}, \dots, n_{{i_k}}\}$ is primitive, let
 $Z_{\{n_{{i_1}}, \dots, n_{{i_k}}\}}$ be the linear subvariety of the affine space ${\bf A}^\rho_k$ defined by $x_{{i_1}} = 0, \dots, x_{{i_k}} = 0$.
 
\noindent  $(\imath \imath \imath)$ Finally, let 
 $$Z_\Sigma =\bigcup_{\{n_{{i_1}}, \dots, n_{{i_k}}\} \rm{~primitive}} Z_{\{n_{{i_1}}, \dots, n_{{i_k}}\}},$$
 called the exceptional set.
 }
 
 \bigskip
 
 The torus ${\bf G}_m^\rho$  over $k$ acts on the affine space ${\bf A}^\rho$ over $k$  diagonaly, by
 $$(\mu_1, \dots, \mu_\rho).(x_1, \dots, x_\rho) := (\mu_1 x_1, \dots, \mu_\rho x_\rho).$$
 The subgroup $G_\Sigma$ of  ${\bf G}_m^\rho$ acts in the same way by restriction, and leaves stable the open set
 ${\bf A}^\rho-Z_\Sigma$. Hence, one can define:

 \medskip

\noindent  {\bf Definition} {\it The toric variety associated to $\Sigma$ is the quotient variety
 $${\bf P}_\Sigma= {\bf A}^\rho-Z_\Sigma /G_\Sigma.$$
 If $(x_1, \dots, x_\rho)$ lies in ${\bf A}^\rho-Z$, we denote by $[x_1, \dots, x_\rho] \in {\bf P}_\Sigma$ its class modulo $G$.
 }

\bigskip

Note that the quotient ${\bf P}_\Sigma$ exists without restriction on the characteristic, thanks to [7]. Since $\Sigma$ is simplicial, a theorem of Cox asserts that this definition coincide with the usual one, given for instance in [8], or in [5].

 \bigskip

\noindent  {\bf II-2. Sub-varieties of simplicial toric varieties.}  Let $\Sigma$ be a simplicial fan in ${\bf R}^d$, with $1$-dimensional cones $n_1, \dots, n_\rho$. Suppose that the group $G_{\Sigma}$ is isomorphic to
  ${\bf G}_m^r$
   with action on ${\bf A}^\rho$ given by
 $$(\mu_1, \dots, \mu_r).(x_1, \dots, x_\rho) = (\prod_{j=1}^r \mu_j^{a_{1, j}}x_j, \dots, \prod_{j=1}^r \mu_j^{a_{\rho, j}}x_j)$$
 for some $a_{i, j} \in {\bf N}$, with $1 \leq i \leq \rho$, $1 \leq j \leq r$. 
 
 Consider the polynomial ring $k[X_1, \dots, X_\rho]$ with one variable $X_i$ for each $n_i$ for $1\leq i \leq \rho$. The ring $ k[X_1, \dots, X_\rho]$ is endowed with a multigraduation
 $$\deg X_i = (a_{i, 1}, \dots, a_{i, r}) \in {\bf N}^r.$$

 Let $P \in k[X_1, \dots, X_\rho]$.  The following remark is crucial:
 
 \medskip
  {\it The condition $P(x_1, \dots, x_\rho)=0$ is independant of the choice of the representative  $(x_1, \dots, x_\rho)$ in ${\bf A}^\rho-Z_\Sigma$  of the class $[x_1 : \dots : x_\rho]$ in ${\bf P}_\Sigma$ if, and only if, $P$ is {\it homogeneous} for the associated multigraduation.
 }
 
 \bigskip
 
 Hence, we can associate to each (multi)-homogeneous ideal $I$ of $k[X_1, \dots, X_\rho]$ a subvariety of ${\bf P}_\Sigma$

 $$V(I) := \{[x_1, \dots, x_\rho] \in {\bf P}_\Sigma, \rm{~such~that~} P(x_1, \dots, x_\rho)=0 \rm{~for~all~} P\in I\}.$$
 
 \bigskip

\noindent  {\bf II-3. Examples.} In the following, we denote by $(e_1, \dots, e_d)$ the canonical basis of
  ${\bf R}^d$.
  
  \medskip
 
\noindent  {\bf II-3-1. The projective space over $k$.} Let $n_1=e_1, \dots, n_d=e_d$, and $n_0=-(e_1+\dots + e_d)$.
 For each free familly $I \subset \{n_0, n_1, \dots, n_d\}$, let $\sigma_I$ be the convex polyedral cone generated by $\{e_i, i \in I\}$. The collection of these $\sigma_I$ form a fan $\Sigma$. Here,
 one as $\rho=d+1$. An elementary calculation shows that
 $G_\Sigma = \{(\mu, \mu, \dots, \mu) \in {\bf G}_m^{d+1}\}$,
 so that $G_\Sigma \simeq {\bf G}_m$, whose action on 
 ${\bf A}^\rho={\bf A}^{d+1}$ is given by
 $\mu.(x_0, x_1, \dots, x_d) = (\mu x_0, \dots, \mu x_d)$. Finally, the only primitive set is
 $\{n_0, n_1, \dots, n_d\}$, so that the exceptional set is $Z_\Sigma=\{(0, \dots, 0)\}$. Hence for this fan, one
 have  
 $${\bf P}_\Sigma = {\bf A}^{d+1}-\{(0, \dots, 0)\}/{\bf G}_m = {\bf P}^d.$$
 As is well known, a polynomial $P(X_0, X_1, \dots, X_d)$ defines a subset of ${\bf P}^d$ if it is homogeneous for the standard graduation.

 \medskip
 
 For positive integers $a_0, a_1, \dots, a_d$, the weighted projective space
 ${\bf P}(a_0, a_1, \dots, a_d)$ is defined in the same way with
 $n_1={1 \over a_1}e_1, \dots, n_d={1 \over a_d}e_d$, and $n_0=-{1 \over a_0}(e_1+\dots + e_d)$,
 giving a fan $\Sigma_a$.
 Here, $G_{\Sigma_a} \simeq {\bf G}_m$ with action on ${\bf A}^{d+1}$ given by
 $\mu.(x_0, \dots, x_d) = (\mu^{a_0}x_0, \dots, \mu^{a_d}x_d)$. Hence,
  a polynomial $P(X_0, X_1, \dots, X_d)$ defines a subset of ${\bf P}(a_0, a_1,$ $\dots, a_d)$ if it is homogeneous for the graduation where $\deg X_i = a_i \in {\bf N}$.
  
  \bigskip
  
\noindent  {\bf Example: Varieties of type III.} Let $n \geq 2$, $d \geq 2$ and $d' \geq 2$ three integers such that $d'$ divides $d$. Fix an hypersurface $X = X_d^{(n-1)}$ of degree $d$ in ${\bf P}^n$ whose equation is $P(x_0, \dots, x_n)=0$. As stated in the introduction, there exists (only one if $k$ is algebraically closed) an $n$-dimensional unramified covering $Y=Y^{(n)}_{X, d'}$ of ${\bf P}^n$, whose ramification divisor is  $X$. 
One can write the equation of $Y$ in the weighted projective space of dimension  $n+1$:
 $$Y = \{[x_0, \dots, x_n, y] \in {\bf P}^{(n+1)}(1, \dots, 1, d/d') ; y^{d'}=P(x_0, \dots, x_n)\}.$$

 \bigskip

\noindent  {\bf II-3-2. Blowing up of ${\bf P}^2$ along a point.} To construct the blowing up of ${\bf P}^2$
 along $[0, 0, 1]$, we start with the fan $\Sigma$ defining ${\bf P}^2$
 in example II-3-1, which have three maximal cones $\sigma_{\{n_1, n_2\}}$,
 $\sigma_{\{n_1, n_0\}}$ and $\sigma_{\{n_2, n_0\}}$. For simplicity, we denote by $[x, y, z]$ the class in ${\bf P}^2$ of the point
 $(x, y, z) \in {\bf A}^3-\{(0,0,0)\}$. 
 We define a new fan by keeping the last two maximal cones (and their faces), add another vector $n_{3} = e_1 + e_2$, an split the maximal cone
 $\sigma_{\{n_1, n_2\}}$ into $2$
 cones
 $\sigma_{\{n_1, n_3\}}$ and $\sigma_{\{n_2, n_3\}}$.
 We obtain a new fan $\widetilde \Sigma$. Now, $\rho = 4$,
 and
 it is easily seen that $(\mu_1, \mu_2, \mu_3, \mu_4) \in G_{\widetilde{\Sigma}}$ if and only if
 $\mu_1 = \mu_2$, and $\mu_3=\mu_1\mu_4$.
 If we denote $\mu=\mu_1=\mu_2$ and $\nu = \mu_4$, we have $\mu_3 = \mu \nu$.
 Hence, $G_{\widetilde{\Sigma}} \simeq {\bf G}_m^2$, whose action on $(x, y, z, v) \in {\bf A}^4$ is
 given by $(\mu, \nu).(x, y, z, v) = (\mu x, \mu y, \mu \nu z, \nu t)$.
 Moreover, the primitive sets are $\{n_1, n_2\}$ and $\{n_0, n_3\}$, hence the exceptional set
 is $Z_{\widetilde \Sigma} = \{(0, 0)\} \times {\bf A}^2 \cup  {\bf A}^2 \times \{(0, 0)\}$.
 
 Then ${\bf P}_{\widetilde{\Sigma}}$ is the blowing up of ${\bf P}^2$ along $[0, 0, 1]$.
 Indeed, let
 $\pi : {\bf P}_{\widetilde{\Sigma}} \rightarrow{\bf P}^2$, given by
 $\pi(x, y, z, v) := (vx, vy, z)$. One observe that this map is well defined, because changing
 $(x, y, z, t)$ by $(\mu, \nu).(x, y, z, t, v)$ changes the image by $(\mu \nu x, \mu \nu y, \mu \nu z)$.
 Moreover, this is an isomophism outside $[0, 0, 1]$, and the inverse image of $[0, 0, 1]$ in
 ${\bf A}^4 - Z_{\widetilde \Sigma}/G_{\widetilde{\Sigma}}$ is isomorphic to ${\bf P}^1$.
 
 Finally, a polynomial $P(X, Y, Z, V)$ defines a subset of ${\bf P}_{\widetilde{\Sigma}}$
  if it is homogeneous for the graduation where $\deg X = \deg Y = (1, 0)$,
  $\deg Z = (1, 1)$ and $\deg V = (0, 1) •n {\bf N}^2$.

 \bigskip
 
\noindent  {\bf II-3-3. Blowing up of  ${\bf P}^4$ along a line.} Let
 $\Sigma$ be the fan in ${\bf R}^4$ defining ${\bf P}^4$, which contains $5$ cones of dimension one $\sigma_{n_i}$ and $5$ maximal simplicial  cones
 $\sigma_{\{n_0, \dots, \widehat{n_i}, \dots, n_4\}}$ for $0 \leq i \leq 4$. 
 Let $L$ be the line $x_1=x_2=x_3=0$.
 We will construct the blowing up of ${\bf P}^4$ along $L$ as a toric variety.
 We add the vector
 $n_5 = e_1 + e_2+e_3$, and we split both maximal simplicial cones $\sigma_{\{n_1, n_2, n_3, n_4\}}$ and
 $\sigma_{\{n_1, n_2, n_3, n_0\}}$ into $3$ cones by replacing each $n_i$, for $i=1, 2$ and $3$, by $n_5$.
 We obtain a new fan $\widetilde{\Sigma}_L$ with $6$ cones of dimension $1$, namely
 $n_i$ for $0 \leq i \leq 5$, and $3 + 3 + 3 = 9$ maximal simplicial cones of dimension $4$.
 
 It is an easy mater to see that $(\mu_i)_{0 \leq i \leq 5} \in G_{\widetilde{\Sigma}_L}$ if, and only if,
 $\mu_1=\mu_2=\mu_3 = \lambda$, $\mu_5 = \mu$ and $\mu_0=\mu_4 = \lambda \mu$.
 Hence, $G_{\widetilde{\Sigma}_L} \simeq {\bf G}_m^2$, with action on ${\bf A}^6$ given by
 
$$(\lambda, \mu).(x_0, \dots, x_5) = (\lambda \mu x_0, \lambda x_1, \lambda x_2,
\lambda x_3, \lambda \mu x_4, \mu x_5).$$
Moreover, the exceptional sets are $\{n_1, n_2, n_3\}$ and $\{n_0, n_4, n_5\}$,
hence
$$Z_{\widetilde{\Sigma}_L} = \{(0, x_1, x_2, x_3, 0, 0)\} \cup \{(x_0, 0, 0, 0, x_4, x_5)\}.$$

The blowing up morphism 
$\pi = {\bf P}_{\widetilde{\Sigma}_L} \rightarrow {\bf P}^4$ is given by

$$\pi((x_0, \dots, x_5)) = (x_0, x_5x_1, x_5x_2, x_5x_3, x_4). \leqno (II.3.3.1)$$
One observe that this map is well defined, is an isomorphism outside $L$, and that the inverse image
of $L$ in ${\bf A}^6-Z_{\widetilde{\Sigma}_L}/G_{\widetilde{\Sigma}_L}$ is isomorphic to ${\bf P}^2\times{{\bf P}^1}$.

 Finally a polynomial $P(X_0, \dots, X_5)$ defines a subset of ${\bf P}_{\widetilde{\Sigma}_L}$
  if it is homogeneous for the graduation where $\deg X_1 = \deg X_2 = \deg X_3 = (1, 0)$,
  $\deg X_5 = (0, 1)$ and $\deg X_0 = \deg X_4 = (1, 1) \in {\bf N}^2$.

 \bigskip

 \noindent {\bf Example : Varieties of type IV.} As stated in the introduction, we will consider for the sake of simplicity only the case of a $3$-fold 
 $X = X_5^{(3)} \subset {\bf P}^4$ containing a line $L$ with multiplicity $3$.
 Denoting by $[x_0, x_1, x_2, x_3, x_4]$ the homogeneous coordinates in ${\bf P}^4$, one can assume, up to a linear change of variables, that $L$ is the line $x_1=x_2=x_3=0$, that is the $x_4$ axis.

 The fact that $L$ have multiplicity $3$ in $X$ means that the affine equation of $X$ in the open set
 $x_0=1$ is of the form
 
 $$P_3(x_1, x_2, x_3) + P_4(x_1, x_2, x_3, x_4) + P_5(x_1, x_2, x_3, x_4) = 0,$$
 where $P_i$ is an homogeneous polynomial of degree $i$ for the standard graduation. Now, the assumption that $L$ is contained in $X$ implies that $x_4$ divides both $P_4(x_1, x_2, x_3, x_4)$ and $P_5(x_1, x_2, x_3, x_4)$.
 Hence, $X$ have homogeneous equation (for the standard graduation)
 
 $$P(x_0, x_1, x_2, x_3, x_4)
 = x_0^2P_3(x_1, x_2, x_3)+x_0x_4 Q_3(x_1, x_2, x_3)+x_4Q_4(x_1, x_2, x_3)
 =0$$
 where $P_3$ and $Q_3$ are homogeneous of degree $3$, and $Q_4$ is homogeneous of degree $4$ for the standard graduation.
 
 \medskip
 
Now, if $\pi$ is the blowing up map given by formula $(II.3.3.1)$, we have:

$$\eqalign{
&[x_0, x_1, x_2, x_3, x_4, x_5] \in \pi^{-1}(X) \cr
& \iff P(x_0, x_1x_5, x_2x_5, x_3x_5, x_4)=0 \cr
 &\iff x_5^3\left\{x_0^2P_3(x_1, x_2, x_3) + x_0x_4Q_3(x_1, x_2, x_3)+x_4x_5Q_4(x_1, x_2, x_3)\right\}=0. \cr
}$$

\noindent By removing three times the strict transform $x_5=0$ of $L$ from this set, we obtain that the strict transform of $X$ in
${\bf P}_{{\widetilde \Sigma}_L}$ has equation
$$\widetilde X_5^{(3)} : x_0^2P_3(x_1, x_2, x_3) + x_0x_4Q_3(x_1, x_2, x_3)+x_4x_5Q_4(x_1, x_2, x_3)=0.
\leqno (II.3.3.2)$$
One observe that this is an homogeneous equation of bi-degree $(5, 2) \in {\bf N}^2$.

 \bigskip

 \noindent {\bf III. Arithmetic over finite fields.} In this chapter, $k$ is the finite field with $q$ elements.
 
  \medskip
 
\noindent  {\bf  III-1. Multigraduate Chevalley-Warning theorem.} We will begin by showing that the corollary of Esnault theorem is easy to prove for varieties of type II, III and IV, and that this proof works also in the singular case.

  \medskip
 
\noindent  {\bf Theorem 1.}  {Ê\it Let ${\bf F}_q[X_1, \dots, X_\rho]$ be the polynomial algebra graduated by ${\bf Z}^r$. Suppose that $\deg X_i \in {\bf N}^r$ for $1\leq i \leq \rho$. Let
 $$(a_1, \dots, a_r)= \deg X_1 + \dots + \deg X_\rho.$$
 Let $P \in {\bf F}_q[X_1, \dots, X_\rho]$, of multidegree $(d_1, \dots, d_r)$. 
 If there exists $j$ such that $d_j < a_j$, then the number of solutions $(x_1, \dots, x_\rho) \in k^\rho$ of the equation $P(x_1, \dots, x_\rho)=0$ is a multiple of $p$. }
 
 \bigskip
 
 {\it Proof.}
 Let $(a_{i, j})_{1 \leq j \leq r} = \deg X_i$ for $1 \leq i \leq \rho$, so that $a_j = \sum_{i=1}^\rho a_{i, j}$.
 We have, with $\alpha = (\alpha_1, \dots, \alpha_\rho)$:
 $$P(x_1, \dots, x_\rho)^{q-1}=
 \sum_{\alpha_1a_{1, j} + \dots + \alpha_\rho a_{ \rho, j} \leq d_j(q-1), 1 \leq j \leq r}
 a_{\alpha} x^{\alpha}$$
 where $x^{\alpha} = x_1^{\alpha_{1}} \times \dots \times x_\rho^{\alpha_{\rho}}$. 
 Hence, if $N = \sharp \{(x_1, \dots, x_\rho) \in {\bf A}^\rho ({\bf F}_q) ; P(x_1, \dots, x_\rho) = 0\}$, we have
 
 $$\eqalign{
 N &= \sum_{x \in {\bf F}_q\rho} (1-P(x)^{q-1})\cr
 &\equiv - \sum_{x \in {\bf F}_q\rho} P(x)^{q-1} ~~~~~~~~~~~~~~~~~ ~~~~~~~~~~~~~~~~~~~~~~~~~~(mod~p)\cr
& \equiv -\sum_{\alpha_1a_{1, j} + \dots + \alpha_\rho a_{\rho, j} \leq d_j(q-1) ; 1 \leq j \leq r}
 a_{\alpha}\prod_{i=1}^\rho (\sum_{x_i \in {\bf F}_q} x_i^{\alpha_i}) ~~(mod~p).\cr
 }$$
 
\noindent But as is well known, for $\alpha \in {\bf N}$, we have
 $\sum_{x \in {\bf F}_q} x^{\alpha} = q-1$ if $\alpha \in {\bf N}^*(q-1)$, and is zero else. Hence, if
 $\alpha_1a_{1, j_0} + \dots + \alpha_\rho a_{ \rho, j_0}
 \leq d_{j_0} (q-1) < a_{j_0}(q-1)
 \leq (q-1)a_{1, j_0} + \dots + (q-1) a_{ \rho, j_0}$ for some $j_0$, we have
 that $0 \leq a_{i, j_0} < q-1$ for at least one index $i \in \{1, \dots, \rho\}$, so that each product in the above formula for $N$ vanishes, hence $N \equiv 0 ~~(mod~p)$.

 \bigskip
 
\noindent  {\bf III-2. Multigraduate Ax theorem.} With more efforts, but with a similar proof than in the classical case, we can prove an Ax theorem for varieties of type II, III et IV,

 \medskip
 
\noindent {\bf Theorem 2.}  {Ê\it Let ${\bf F}_q[X_1, \dots, X_\rho]$ be the polynomial algebra, graduated by ${\bf Z}^r$. Suppose that $\deg X_i \in {\bf N}^r$ for $1\leq i \leq \rho$. Let
 $$(a_1, \dots, a_r)= \deg X_1 + \dots + \deg X_\rho.$$
 Let $P \in {\bf F}_q[X_1, \dots, X_\rho]$, of multidegree $(d_1, \dots, d_\rho)$ and
 $$N = \sharp \{(x_1, \dots, x_\rho) \in {\bf A}^\rho({\bf F}_q) ; P(x_1, \dots x_\rho) = 0\}.$$ 
 Then
 $N\equiv 0 ~(mod~q^\mu)$,
 where
 $\mu = \max_{1\leq i\leq r} \Big\lceil {a_i-d_i \over d_i}Ê\Big\rceil$.
  }

 \bigskip
 
 {\it Proof.} We will prove only for simplicity that $N \equiv 0~(mod~q^{\mu_1})$ for $\mu_1 = \Big\lceil {a_1-d_1 \over d_1}Ê\Big\rceil$. One may assume that $d_1 < a_1$. 
We will adapt Ax's original proof given in [1]. We will begin by collecting in lemmas 3, 4, 5 and 6 all results from [1] whose proofs don't use any consideration on the graduation. We will then conclude for the multigraduations case.

 Let $(a_{i, j})_{1 \leq j \leq r} = \deg X_i$ for $1 \leq i \leq \rho$. First of all, we recall (and modify slightly) Ax's notations. If

$$W = \Big\{w =(w_1, \dots, w_rho) \in {\bf N}^\rho
\mid weight(w) : = a_{1, 1}w_1 + \dots a_{\rho, 1}w_\rho \leq d_1\Big\},$$

\noindent then one can write $P(X) = \sum_{w \in W} a(w) X^{w}$ for some
$a(w) \in {\bf F}_q$.
Let $T$ be the Teichm\"uller set of representatives of ${\bf F}_q$
in the unramified extension $K$ of ${\bf Q}_p$ whose residue field is ${\bf F}_q$.
Let $\zeta$ be a primitive $p$-th root of unity in $\overline {\bf Q}_p$.
By Lagrange interpolation theorem, there is an unique polynomial

$$C(U) = \sum_{m=0}^{q-1} c(m)U^m \in K(\zeta)[U],$$

\noindent such that $C(t) = \zeta^{Tr_{{\bf F}_q/{\bf F}_p} (t)}$ for any $t \in T$.

\medskip

\noindent {\bf Lemma 3 (Ax).} {\it We have

$$c(m) \equiv~0~mod~(1-\zeta)^{\sigma(m)}$$
for any $0 \leq m \leq q-1$, where $\sigma(m)$ is the $p$-weight of $m$.}

\bigskip

\noindent {\it Proof of lemma 3.} This is relation (4) in [1]. It follows from Stickelberger theorem.

\bigskip

 Let $M$ be the set of functions $m$
from $W$ to $\{0, \dots, q-1\}$. For $m \in M$, let

$$e(m) = \sum_{w \in W} m_i(w)~ w \in {\bf N}^\rho,$$
and

$$e(m)' = \sum_{w \in W} m_i(w)  \in {\bf N}.$$

\medskip

\noindent {\bf Lemma 4 (Ax).} {\it
$$\eqalign{
q N &=
\sum_{m \in M} \left\{\prod_{w \in W} A(w)^{m_i(w)}\right\}\cr
&~~~~~~\times \left\{\prod_{w \in W} C(m_i(w))\right\}\cr
&~~~~~~\times \left\{\Bigl(\sum_{t \in T^\rho} t^{e(m)}\Bigr).
\Bigl(\sum_{t' \in T} t'^{e'(m)}\Bigr)\right\},
}$$
where $A(w)$ is the Teichm\"uller representative of $a(w)$ in $T$.}

\bigskip

\noindent {\it Proof of lemma 4.} This is relation (5') in [1]. The proof involves the existence
of an additive character $\beta$ from ${\bf F}_q$ with values in $K$.

\bigskip

\noindent {\bf Lemma 5 (Ax).} {\it Let $m \in M$.

\noindent $(\imath)$ If $e(m)$ and $e(m)'$ are not both multiples of  $q-1$ in ${\bf N}^\rho$
and ${\bf N}$ respectively, then

$$\sum_{t \in T^\rho} t^{e(m)} .\sum_{t' \in T}
t'^{e'(m)} = 0.$$

\noindent $(\imath \imath)$ Suppose that $e(m)$ and $e'(m)$ are multiples of $q-1$.
Denote by $s_1 \in \{0, 1, \dots, \rho\}$ the number of non-zero entries of $e(m)$,
and by $s_2 \in \{0, 1\}$ those of $e'(m)$. Then

$$\sum_{t \in T^\rho} t^{e(m)} .\sum_{t' \in T}
t'^{e'(m)} = (q-1)^{s_1 + s_2}q^{\rho+1-s_1-s_2}.$$
}

\bigskip

\noindent {\it Proof of lemma 5.} This is relations (6), (6') and (6'') in [1]. 
\bigskip

The statement of lemma 6 requiers the following notations. If $n$ is any integer, and if $p$ is the characteristic of ${\bf F}_q$, we write the decomposition of any integer $n \in {\bf N}$ in bases $p$ as
$$n=\sum_{i \in {\bf N}} n_ip^i,$$
where $0 \leq n_i \leq p-1$. Then, we define $\sigma(n) = \sum_{i \in {\bf N}} n_i$. Let $f \in {\bf N}^*$ be such that
$q=p^f$.
\medskip

\noindent {\bf Lemma 6 (Ax).} {\it Let $m \in M$, and $j_1, \dots, j_{s_1}$ be the indices for which $e(m)_j \neq 0$. We have
$$f(p-1) \Bigg\lceil{\sum_{k=1}^{s_1} a_{j_k, 1} -d_1\over d_1}\Bigg\rceil
\leq \sum_{w \in W} \sigma(m(w)).$$
}

\medskip

\noindent {\it Proof of lemma 6.} The proof follows the same pattern than formula preceding formula (8) one in [1].

\bigskip

\quad We are now ready to prove theorem 2. From lemmas 4 and 5, the functions $m$  whose
contribution in the number $N$ are non-trivial are
those for which $e(m)$ and
$e'(m)$ are both multiples of $q-1$ by non-zero integer-valued vectors. Let
$m$ be such a function and $j_1, \dots, j_{s_1}$ be the indices for which $e(m)_j \neq 0$. Now, lemmas 3, lemma 6 and the fact that $p$ divides $(1-\zeta)^{p-1}$ imply that the exponant of the greater power of $q$ dividing $\prod_{w \in W} c(m(w))$  is at least
$$\Bigg\lceil{\sum_{k=1}^{s_1} a_{j_k, 1} -d_1\over d_1}\Bigg\rceil.$$
Together with lemma 5, lemma 4 implies that
$$q^r \vert qN$$
with 
$$r = \min_{0 \leq s_1 \leq \rho ; \{a_{j_1, 1}, \dots a_{j_{s_1}, 1}\} \subset \{a_{1, 1}, \dots a_{\rho, 1}\}} r(j_1, \dots, j_{s_1}),$$
 where
$$r(j_1, \dots, j_{s_1}) = \Bigg\lceil{\sum_{k=1}^{s_1} a_{j_k, 1} -d_1\over d_1}\Bigg\rceil + \rho - s_1.$$
Now, for any $\{a_{k_1, 1}, \dots a_{k_h, 1}\}$ disjoint to $\{a_{j_1, 1}, \dots a_{j_{s_1}, 1}\}$, we have
$$h \geq
\Bigg\lceil{\sum_{i=1}^{s_1} a_{j_i, 1} + \sum_{i=1}^{k} a_{k_i, 1}-d_1\over d_1}\Bigg\rceil
- \Bigg\lceil{\sum_{i=1}^{s_1} a_{j_i, 1} -d_1\over d_1}\Bigg\rceil.$$
Subsituting $h = \rho-s_1$ and
$\{a_{k_1, 1}, \dots a_{k_h, 1}\}=\{a_{1, 1}, \dots a_{\rho, 1}\}-\{a_{j_1, 1}, \dots a_{j_{s_1}, 1}\}$, we obtain
that $r(j_1, \dots, j_{s_1}) \geq  \Bigg\lceil{\sum_{i=1}^{\rho} a_{i, 1} -d_1\over d_1}\Bigg\rceil = \lceil{a_1-d_1 \over d_1}\rceil$, which was to be proved.

 \bigskip
 
 \noindent {\bf Remark.}  It is likely that a congruence relation {\it \`a la} Katz can be given for several polynomial equations, adapting Wan's proof given in [13] (see [11] for non-standard graduation by ${\bf N}$).

 \bigskip
 
\noindent  {\bf III-3. Back to Esnault theorem for RC-quintic $3$ folds.} As we saw in section II.3.3, formula
 II.3.3.2, the equation $P(x_0, \dots, x_5) = 0$ of a RC-quintic $3$-fold is homogeneous of degree $(d_1, d_2) = (5, 2)$ in the
 toric variety ${\bf P}_{\widetilde{\Sigma}_L}$ of section II.3.3. Here, $r=2$ and $\rho = 6$, and
 $(a_1, a_2) = \deg X_0 + \dots + \deg X_5 = (1, 1)+ 3 \times (1, 0) + (1, 1) + (0, 1) = (5, 3)$.
 Since $d_2 = 2 < a_2=3$, we have $\mu = \lceil {3-2 \over 2}\rceil = 1$, and multigraduated Ax theorems asserts that the number of affine points
 $N_{aff} = \sharp \{(x_0, \dots, x_5) \in {\bf A}^6({\bf F}_q) ; P(x_0, \dots, x_5) = 0\}$ is congruent to zero $(mod~q^\mu = q)$.
 
 Now, we whant to study the number of rational points of the complete variety
 $\widetilde X_5^{(3)} \in {\bf P}_{\widetilde{\Sigma}_L} = {\bf A}^6-Z_{\widetilde{\Sigma}_L}/{\bf G}_{\widetilde{\Sigma}_L}$. The exceptional set $Z_{\widetilde{\Sigma}_L}$ is given in section $III.3.3$.
 It is easily seen that it is contained in the affine zero set $P(x_0, \dots, x_5) = 0$ where $P$ is given by formula II.3.3.2. Moreover, ${\bf G}_{\widetilde{\Sigma}_L} \simeq {\bf G}_m^2$, hence

 $$\eqalign{
 \sharp \widetilde X_5^{(3)}({\bf F}_q)
&\equiv {N_{aff}-\sharp Z_{\widetilde{\Sigma}_L}\over (q-1)^2} \cr
&\equiv {0-(q^3+q^3-1) \over (q-1)^2} \cr
 &\equiv 1~(mod~q),\cr
}$$
as stated by Esnault theorem in this case.

\bigskip

 \noindent {\bf IV. Arithmetic over $k(T)$, $k$ algebraically closed.} Always for the sake of simplicity, we will consider only the case of a RC quintic $3$-fold $ \widetilde X_5^{(3)}$ of type IV over $k(T)$, where $k$ is an algebraically closed field. Of course, the following proof can be extended to many other subvarieties of ``low degree" in some toric varieties. It is an extension in this toric case of the classical Tsen's proof.
  
 \medskip
 
\noindent  {\bf Theorem 3.} {\it Let $\widetilde X_5^{(3)} \in {\bf P}_{\widetilde{\Sigma}_L}$ be a RC-quintic $3$-fold, with equation $P(x_0, \dots, x_5)=0$ of the form $II.3.3.2$ with coefficients in $k(T)$, where $k$ is an algebraically closed field. Then $\widetilde X_5^{(3)}(k(T)) \neq \emptyset$.
 }
 
 \bigskip
 
\noindent  {\it Proof of theorem 3.} To begin with, recall the general equation II.3.3.2 of $\widetilde X_5^{(3)} \in {\bf P}_{\widetilde{\Sigma}_L}$:
 $$x_0^2P_3(x_1, x_2, x_3) + x_0x_4Q_3(x_1, x_2, x_3)+x_4x_5Q_4(x_1, x_2, x_3)=0,$$
 where all polynomials have coefficients in $k(T)$, and are homogeneous for the standard graduation, 
 of degree given by their indice. Recall also that $k[X_0, \dots, X_5]$ is graduated by
 $\deg X_1 = \deg X_2 = \deg X_3 = (1, 0)$,
 $\deg X_5 = (0, 1)$ and $\deg X_0 = \deg X_4 = (1, 1) \in {\bf N}^2$.
 
As in original Tsen proof, we extend each variable $(x_0, \dots, x_5)$
into polynomial variables 
$$x_i(T) = \sum_{j=0}^s x_{i, j}T^j,$$
for $0 \leq i \leq 5$, with $x_{i, j} \in k$.
We have to prove that there exists
$(x_0(T), \dots , x_5(T)) \in {\bf A}_{k(T)}^6 - Z_{\widetilde{\Sigma}_L}(k(T))$, such that

$$x_0(T)^2P_3(x_1(T), x_2(T), x_3(T)) + x_0(T)x_4(T)Q_3(x_1(T), x_2(T), x_3(T))$$
$$+x_4(T)x_5(T)Q_4(x_1(T), x_2(T), x_3(T))=0.$$
This expression have the form
 
 $$\sum_{n=0}^{5s+c}F_n(x_{i, j}, 1 \leq i \leq 5 ; 1 \leq j \leq s) =0,$$
 where $F_n(X_{i, j})$ is an homogeneous polynomial of degree $(5, 2)$ for the graduation
 $\deg x_{i, j} = \deg X_i$. This equation is equivalent to the system with $5s + c + 1$ equations and $6s+6$
 unknowns:

 $$\left \{
\eqalign{ 
&F_n(x_{i, j}) = 0 \cr
&0 \leq n \leq 5s+c,\cr
}\right .$$
where we have to avoid the exceptional set
$$Z_{\widetilde{\Sigma}_L}(k(T)) = \{x_1(T)=x_2(T)=x_3(T)=0\}\cup \{x_0(T)=x_4(T)=x_5(T)=0\}.$$

Now, the rational chow ring $A_{\bf Q}^*({\bf P}_{\widetilde{\Sigma_L}})$ is isomorphic to (this follows from [8], chapter 5):

$$A_{\bf Q}^*({\bf P}_{\widetilde{\Sigma}_L})
= {\bf Q}[x_0, \dots, x_5]/
<x_1=x_2=x_3,
x_0=x_4=x_1+x_5,
x_1x_2x_3=0,
x_0x_4x_5=0>,$$
so that denoting $x = x_1=x_2=x_3$, $v=x_5$ and $u=x_0=x_4$, we have
$$A_{\bf Q}^*({\bf P}_{\widetilde{\Sigma}_L})
= {\bf Q}[x, u, v]/
<u=x+v, x^3=0, u^2v=0>.$$

Since the zero set of an homogeneous polynomial of degree $(5, 2)$ is equivalent in
$A_{\bf Q}^*({\bf P}_{\widetilde{\Sigma}_L})$ to $3x+2u=4x+u+v$, we have only to prove that
$(3x+2u)^{5s+c+1} \neq 0$ in the new ring
$$A_{\bf Q}^s := {\bf Q}[x, u, v]/
<u=x+v, x^{3s+3}=0, u^{2s+2}v^{s+1}=0>.$$
In fact, we have in $A_{\bf Q}^s$, for $s$ large enought:

$$\eqalign{
(3x+2u)^{5s+c+1}\times v^{s-c+3}
&=(u+v+4x)^{5s+c-1}\times v^{s-c+3} \cr
&= \gamma x^{3s+2}\times (u+v)^{2s+c-2}\times v^{s-c+3} \cr
&= \gamma x^{3s+2}\times u^{2s+2}v^s\cr
& \neq 0,
}$$
where $\gamma$ is a positive integer linear combination of non-zero binomial coefficients. This proves that
$(3x+2u)^{5s+c+1}\times v^{s-c+3}$, hence also $(3x+2u)^{5s+c+1}$, are non-zero in $A_{\bf Q}^s$,
and the proof is complete.

 \bigskip
 
\noindent  {\bf Remark.} To extend this theorem to subvarieties of other toric varieties over $k(T)$, we
 have to extend the variable $x_i$ to the polynomial variable
 $x_{i, 0} + t x_{i, 1} + \dots + x_{i, sa_{i, r}} t^{sa_{i, r}}$ if
 $\deg X_i = (a_{i, j})_{1 \leq j \leq r}$.

\bigskip

\centerline{\bf References.}

\medskip

\noindent [1] J. Ax, {\it Zeroes of polynomial over finite fields}, Amer. Journ. Math. 86 (1964), 255--261.

\medskip

\noindent [2] F. Campana, {\it Connexit\'e rationnelle des vari\'et\'es de Fano}, Ann. Sc. \'Ec. Norm. Sup. 25 (1992), 539--545.

\medskip

 \noindent [3] A. Chambert-Loir, {\it Points rationnels et groupes fondamentaux : applications de la cohomologie $p$-adique}, S\'eminaire Bourbaki 55\`e ann\'ee, 2002-2003, no 914, mars 2003.

\medskip

 \noindent [4] D. Cox, {\it The homogeneous coordinate ring of toric varieties},  Journ. of Alg. Geom. 4 (1995), 17--50.

\medskip 

\noindent [5] V. I. Danilov, {\it The geometry of toric varieties}, Russian Math Surveys 33 : 2 (1978),97--154.
\medskip

 \noindent [6] O. Debarre, {\it Vari\'et\'es rationnellement connexes}, S\'eminaire Bourbaki 54\`e ann\'ee, 2001-2002, no 905, juin 2002.

\medskip 

\noindent [7] J.. Forgaty, F. Kirwan, D. Mumford, {\it Geometric Invariant Theory}, Ergebnisse der maths  34, Springer-Verlag, 1994.

\medskip

\noindent [8] W. Fulton, {\it Introduction to toric varieties}, Princeton University Press, Princeton, New Jersey, 1993.

\medskip

\noindent [9] J. Kollar, {\it Rational curves on algebraic varieties}, Ergebnisse der maths 32, Springer-Verlag, 1996.

\medskip

\noindent [10] J. Kollar, Y. Miyahoka, S. Mori, {Rationally connected varieties}, J. Alg. Geom.
 1 (1992), 429--448.

\medskip

\noindent [11] M. Perret, {\it On the number of points of some varieties over finite fields}, Bull London Math. Soc. 35 (2003), 309--320.

\medskip 

\noindent [12] J.-P. Serre, {\it Rev\^etements ramifi\'es du plan Projectif}, S\'eminaire Bourbaki 1959-60, expos\'e 204.

\medskip

\noindent [13] D. Wan, {\it An elementary proof of a theorem of Katz}, Amer. journ. Maths, 111
(1989), 1--8.

\bigskip

\noindent Marc Perret

\noindent \'Equipe \'Emile Picard,

\noindent Institut de Math\'ematiques de Toulouse, UMR 5219

\noindent Universit\'e de toulouse II

\noindent 5, all\'ees Antonio Machado

\noindent 31 058 Toulouse Cedex

\noindent perret@univ-tlse2.fr

\end